\documentclass[11pt]{article}
\usepackage{amsfonts}
\usepackage{amssymb}
\usepackage{amssymb,amsmath}
\usepackage{graphicx}
\usepackage{flafter}
\makeatletter

\newcommand{\Rmnum}[1]{\expandafter\@slowromancap\romannumeral #1@}
\makeatother

\def\o{\omega}

\def\R{{\bf R}}
\def\E{{\bf E}}
\def\P{{\bf P}}

\def\C{{\bf C}}

\def\p{\partial}

\newcommand{\be}{\begin{equation}}
\newcommand{\ee}{\end{equation}}
\newcommand{\bd}{\begin{displaymath}}
\newcommand{\ed}{\end{displaymath}}
\newcommand{\ba}{\begin{array}{ll}}
\newcommand{\ea}{\end{array}}
\newcommand{\baa}{\begin{eqnarray}}
\newcommand{\eaa}{\end{eqnarray}}
\newcommand{\baaa}{\begin{eqnarray*}}
\newcommand{\eaaa}{\end{eqnarray*}}
\font\sm=cmr10

\def\T{{\cal T}}

\date{}
\title{
Boundary Value Problems for Functionals of Ito
Processes\footnote{{\it Theory of Probability and its Applications} {\bf
36} (1992), No. 3, 459-476. Translated from Russian Journal  {\em Teoriya Veroyatnostei i ee Primeneniya}  {\bf
36} (1991), No. 3, pp. 464-481}
\footnote{Received by the editors September 27, 1988.}}

\author{
N. G. Dokuchaev\\ {\sm  (Translated by V.A.Lebedev)}}
 
\begin{document}
 \vspace{-0.5cm}      \maketitle
\section{Formulation of the problem and main assumptions}

Let us consider a probability space $(\Omega,\mathcal
{F},\textbf{P})$, where $\Omega = \{\omega\}$ is a set of elementary
events, $\mathcal {F}$ is some $\textbf{P}$-complete
$\sigma$-algebra of events, $\textbf{P}$ is a probability measure on
$\mathcal {F}$. We consider a standard $d_0$-dimensional Wiener
process $W(t)=\|\omega_1(t),\dots,\omega_{d_0}(t)\|$ with
independent components. The part of this process
$\|\omega_1(t),\dots,\omega_d(t)\|$, where $d \leq d_0$, is denoted
by $\omega(t)$. The process $\omega(t)$ generates the filtration of
$\textbf{P}$-complete $\sigma$-algebras
$\mathcal{F}_t=\overline{\sigma[\omega(s), s \leq t]}\subset
\mathcal{F}$ in the usual way.
\par

We consider a $n$-vector It\^o stochastic differential equation
\begin{align*}
    dy^{x,s}(t,\omega)&=f[y^{x,s}(t,\omega),t,\omega]dt +
    \beta[y^{x,s}(t,\omega),t,\omega]dW(t), \tag{1.1}\\
    y^{x,s}(s,\omega)&=x, \tag{1.2}
\end{align*}
where $0\leq s \leq t \leq T$, $x\in R^n$, and the number $T>0$. The
function $f(x,t,\omega):R^n \times R^{+} \times \Omega \rightarrow
R^n$,$\beta(x,t,\omega):R^n\times R^+\times \rightarrow R^{n\times
d}$ are progressively measurable with respect to the filtration of
$\sigma$-algebras $\textit{F}_t$ for any $x\in R^n$. These functions
are measurable, bounded, satisfy the global Lipschitz condition in
$x$ uniformly in $t,w,$ and are continuous in $x,t$ for any $w$. By
a solution of $(1.1)$, $(1.2)$ we shall mean a "strong" solution.

\par
    Let a region $D\subset \textbf{R}^n$ be given, and let either
    $D=\textbf{R}^n$ or the region $D$ be simply connected, bounded,
    and have a $C^2$-smooth boundary. Let us consider the cylinder
    $Q=D\times(0,T)$, and , for each $(x,s)\in \overline{Q}$, the random
    variable $\tau ^{x,s}(\omega)= T \bigwedge
    \inf\{t:y^{x,s}(t,\omega)\notin
    \overline{D}\}$, that is, the first exit time from the set $\overline{Q}=Q\cup \partial Q$
    for the vector$[y^{x,s}(t,\omega),t]$. If $D=\textbf{R}^n$, then
    $\tau ^{x,s}(\omega)\equiv T$.

\par
    This paper is devoted to the study of functionals of the form
\begin{equation}\tag{1.3}
    v(x,s,\omega)= E\left\{ \int_s^{\tau^{x,s}(\omega)}\varphi \left[ y^{x,s}(t,\omega),t, \omega \right] dt | \mathcal {F}_s\right\}
\end{equation}

    Here the functions $\varphi(x,t,\omega:\textbf{R}^n\times \textbf{R}^+\times \Omega\rightarrow
    \textbf{R})$ are progressively measurable with respect to the
    filtration $\mathcal{F}_t$ for any $x\in \textbf{R}^n$;
    $E\{\cdot|\mathcal{F}_s\}$ is the conditional expectation.

\par
    For distributions of such functionals of It$\hat{o}$ processes,
    which are not Markov,estimates are given in $[1,
    \text{Chap}.\Rmnum{2}]$.

\par
    The goal of the paper consists in the representation of
    functional $(1.3)$ by solutions of special boundary value
    problems for stochastic partial differential equations
    introduced in $\S2$. In $\S3$ we establish the duality of these
    problems to boundary value problems for It$\hat{o}$ parabolic
    equations which allows us to obtain supplementary information
    about solutions of boundary value problems of both forms(Theorem
    3.2) and Theorem 4.1). Sufficient conditions for a
    representation of a solution of the boundary value problem in
    the form $(1.3)$ are obtained in $\S 2$, sufficient conditions
    for a representation of the functional $(1.3)$ in the form of a
    solution of a boundary value problem are obtained in $\S 5$ (these
    cases are different because the function $\varphi$ does not coincide
    with the free term of the partial equation if the process $y^{x,s}(t,\omega)$
    is not Markov). A certain smoothness of the functionals (1.3) in
    $x,s$ is also established (Theorem 5.1).

\par
    Let us make additional assumptions.

\par
    For $j=1,\dots,d_0$ we denote  by $\beta_j$ the corresponding
    columns of the matrix $\beta$. In the case $d<d_0$ we denote by
    $\tilde{\beta}$ the $n\times (d_0-d)$-matrix
    $\|\beta_{d+1},\dots,\beta_{d_0}\|$. We assume that the
    eigenvalues of the matrices $\beta\beta^T$ and
    $\tilde{\beta}\tilde{\beta}^T$( in the case $d<d_0$) are
    separated from zero uniformly in all arguments.
\par
    Let us fix an integer number $r\geq 0$ and a number $l>0$ such
    that $r<l<r+1$, $r=[l]$. Let, as in $[3,p.7]$, $H^{l,l/2}(\overline{Q})$
    be the same Banach space of functions on $\overline{Q}$ which are
    H$\ddot{o}$lder continuous together with $r$ derivatives in $x$
    and $[r/2]$ derivatives in $t$. We assume that the functions $f(x,t,\omega)$
    and $\beta(x,t,\omega)$ belong componentwise to
    $H^{l,l/2}(\overline{Q})$ for every $\omega \subset \Omega$ and their
    norms in this space are bounded uniformly in $\omega \in
    \Omega$. For $r<2$ the partial derivatives of the components of
    the matrix $\beta(x,t,\omega)$ of second order in $x$ are
    assumed to be uniformly bounded in $x,t,\omega$. In the case
    $r>0$ and $D\neq \textbf{R}^n$ we assume that the boundary $\partial
    D$ belongs to the class $H^{l+2}$(see $[3,p.9]$). For
    $D=\textbf{R}^n$ we have $\partial D=\Phi$ and by $\overline{D}$ and
    $\overline{Q}$ we mean $\textbf{R}^n$ and $\textbf{R}^n$ and $\textbf{R}^n \times
    [0,T]$, respectively.

\par
    Below, $L_2(D),L_2(Q), W_2^m(D), W_2^1(\overline{D}),
    C^m(\overline{D}),C(\overline{Q})$, and so on denote the usual spaces
    $(|[2]-[5]|)$ of real-valued functions on $\overline{D}$ or
    $\overline{Q}$. For a Banach space $\mathcal {X}$ the symbol
    $\|.\|_{\mathcal{X}}$ denotes the norm, for a Hilbert space
    $\mathcal{X}$ the symbol $(.,.)_{\mathcal{X}}$ denotes the
    scalar  product. For a region $G\subset \textbf{R}^m$ the symbol
    $C(\overline{G}\rightarrow \mathcal{X})$ denotes the Banach space of
    continuous bounded functions $u: \overline{G}\rightarrow \mathcal{X}$
    with the usual norm. $C^{m,q}(\overline{Q}\rightarrow \mathcal{X})$
    denotes the set of functions $u(x,t):\overline{Q}\rightarrow \mathcal{X}$
    belonging to $C(\overline{Q}\rightarrow \mathcal{X})$ together with
    the first derivatives in $x$ and $q$ derivatives in $t$.

\par

    Let us consider the positive self-dual unbounded operator $\Lambda: L_2(D)\rightarrow L_2(D)$
    of the form $\Lambda = \sqrt{I-\Delta}$,
    where $I$ is the identity operator and $\Delta$ is
    $n$-dimensional Laplace operator. For $k=0,\pm 1$ we introduce
    the Hilbert spaces $H^k$ with the scalar product $(u,v)_{H^k}=(\Lambda^k u ,\Lambda^k v)_{L_2(D)}$.
    We assume that $H^{-1}$ is the completion of
    $L_2(D)$ in the norm $\|.\|_{H^{-1}}$, $H^0=L_2(D)$,
    $H^1=W^1_2(R^n)$ for $D=R^n$, and $H^1= W_2^1(R^n)$ for $D=R^n$,
    and $\dot{W}^1_2(D)$ for $D\neq R^n$.  The coincidence of
    the corresponding norms for $k=0,1$ can be easily verified (see the description of $H^k$ in
    $[2]$). For $u\in H^1$ and $v\in H^{-1} $ by $(u,v)_{H^0}$ we mean $(\Lambda u,\Lambda^{-1}
    v)_{H^0}$.

\par
    The symbol $\lambda_1$ denotes the Lebesgue measure in $[0,T]$.
    $ \overline{\mathcal {P}}$ (and $\overline{\mathcal{P}}_s $ for a given $s \in
    [0,T]$) denotes the completion in the measure $\lambda_1 \in
    \mathcal{P}$ of the $\sigma$-algebra of subsets of the set $[0,T]\times\Omega$
    generated by stochastic processes which are progressively
    measurable with respect to the filtration $\mathcal{F}_t$( respectively, of the $\sigma$-algebra generated by measurable processes $\xi(t,\omega)$
    for all $t\in[0,T]$ )which are measurable with respect to
    $\mathcal{F}_s$.

\par
    For integer numbers $m\geq 0$, $k=0,\pm1$, we introduce the
    Hilbert spaces
\begin{align*}
\mathcal{L}_2 &=L^2([0,T]\times\Omega, \overline{\mathcal{P}},
\lambda_1
  \times P,R),\\
  X^k &= L^2([0,T]\times \Omega, \overline{\mathcal{P}}, \lambda_1
  \times P, H^k),\\
  \overline{X}^k &= L^2([0,T]\times \Omega, \overline{\mathcal{P}}_T, \lambda_1\times P,
  H^k),\\
  \mathcal{W}^m &=L^2([0,T]\times\Omega, \overline{\mathcal{P}}, \lambda_1 \times P,
  W_2^m(D)),\\
  \overline{\mathcal{W}}^m &= L^2([0,T]\times \Omega, \overline{\mathcal{P}}_T, \lambda_1\times P,
  W_2^m(D)).
\end{align*}

    For $p\geq 1$, $s\in[0,T]$ and the number $l$ fixed above, we
    introduce the Banach spaces
\begin{align*}
    \overline{\mathcal{H}}^l &= L^2(\Omega, \mathcal{F}_T, P,
    H^{l,l/2}(\overline{Q})),\\
    C_p^m &= L^p([0,T]\times \Omega, \overline{\mathcal{P}}, \lambda_1 \times P,
    C^m(\overline{D})),\\
    \overline{C}_p^m(s) &= L^p([0,T]\times \Omega, \overline{\mathcal{P}}_s, \lambda_1 \times P,
    C^m(\overline{D})),\\
    C_0 &= C([0,T]\rightarrow L^2(\Omega, \mathcal{F}_T, P,
    L_2(D))),\\
    \overline{\mathfrak{C}} &= L^2(\Omega, \mathcal{F}_T, P, C(\overline{Q})).
\end{align*}

    For integer numbers $m \geq 0$, $q \geq 0$, the symbol
    $\overline{C}^{m,q}$ denotes the set of functions $u(x,t,\omega)$
    belonging to $\overline{C}$ together with the first $m$ derivatives
    in $x$, and $q$ derivatives in $t$ (the derivatives must exist with probability
    1).

\par
    We assume that $C_0\subset X^0 \subset X^{-1}$, $X^1 \subset \mathcal {W}^1 \subset X^0 = \mathcal
    {W}^0$, $\overline{\mathcal {H}}^l \subset \overline{\mathcal {C}} \subset
    \overline{\mathcal{C}}_2^r(T)$, and so on, meaning the natural dense
    embedding. Moreover, $\mathcal {C}_p^m \subset \overline{\mathcal
    {C}}_p^m(T)$, $X^k \subset \overline{X}^k$, and so on. $\mathcal
    {H}^l$ denotes the set $\overline{\mathcal {H}}^l\cap \mathcal
    {W}^r$, where $r=[l]$.

\par
    We introduce the set $\partial_0Q \subset \partial Q$ and the
    set $\partial_T Q \subset \partial Q$ of the following form:
\[
    \partial_0 Q = \{\partial D \times [0,T]\} \cup \{D \times
    \{0\}\}, \p_T Q = \{\p D \times [0,T]\} \cup \{D \times \{T\}\}
    ;
\]
    in the case $D \neq R^n$, $\p_0 Q = R^n \times \{0\}$, in the
    case $D = R^n $, $\p_T Q = R^n \times \{\T\}$.

\par
    For every $\omega \in \Omega$ we define the differential
    operator
\begin{equation}\tag{1.4}
    A = A(x,t,\omega)=\sum_{i=1}^n f_i(x,t,\omega)\frac{\p}{\p x_i}
    + \frac{1}{2}\sum_{i,j=1}^n b_{ij}(x,t, \omega) \frac {\p^2}{\p x_i \p
    x_j}.
\end{equation}

    Here $f_i$, $x_i$, $b_{ij}$ are components of the vectors $f$,
    $x$, and of the matrix $b= \beta \beta^T$. $A^*(x,t, \omega)$
    will denote the differential operator dual to the operator(1.4)
    (in the Lagrange sense (see [4, p.141])).

\par
    For $g \in \mathcal{H}^l$ we consider the following boundary
    value problem in $Q$:
\begin{align*}
    \frac{\p U}{\p t}(x,t, \omega) + A(x,t,\omega)U(x,t,\omega) =
    -g(x,t,\omega), \tag{1.5}\\
    U(x,t,\omega)|_{(x,t)\in \p_T Q} = 0. \tag{1.6}
\end{align*}

\par
    We introduce the operator $\overline{\mathcal {T}}$, which maps the
    function $g$ to a solution $U =\overline{\mathcal {T}}g $ of the
    boundary value problem (1.5)-(1.6). From [3](see also [2] and
    [4]) it follows that the operator $\overline{\mathcal {T}}: X^0 \rightarrow \mathcal
    {W}^2$, $\overline{\mathcal {T}}: X^{-1}\rightarrow \overline{X}^1$, $\overline{\mathcal {T}}: X^{-1} \rightarrow
    C_0$ are continuous. Moreover, $U=\overline{\mathcal {T}} g \in
    \mathfrak{C}^{r=2,1}$ if $g \in \mathcal {H}^l$.

\section{Representation of solutions of boundary value problems in the form of functionals of Ito processes}
    In the cylinder $Q$ we consider the following boundary value
    problem for a stochastic partial differential equation:
\begin{align*}
    d_t v(x,t,\omega) + [A(x,t,\omega)v(x,t,\omega) +
    g(x,t,\omega)]dt = \mathcal {X}(x,t,\omega)d \omega(t), \tag{2.1}\\
    v(x,t,\omega)|_{(x,t)\in\p _T Q} = 0. \tag{2.2}
\end{align*}

    Here the function $v$ is scalar-valued and values of the
    function $\mathcal {X}$ are row $d$-vectors, $\mathcal {X} = \|\mathcal {X}_1, \dots, \mathcal
    {X}_d\|$. Equation (2.1) in combination with a boundary
    condition at $t=T$ means, in the case $v \in \mathcal{C}_2^2 \cap
    C_0$, $g\in X^0$, $\mathcal {X}_j \in X^0$, that for any $t$ for
    a.e.(almost every) $x$, $\omega$,
\begin{equation} \tag{2.3}
    v(x,t,\omega) = \int_t^T[A(x,\rho,\omega)v(x,\rho,\omega) +
    g(x,\rho,\omega)]d\rho -\int_t^T\mathcal {X}(x,\rho,\omega)
    d\omega(\rho).
\end{equation}

\par
    The stochastic integral with respect to $d\omega_j(\rho)$ of a
    square-summable progressively measurable with respect to the
    filtration $\mathcal{F}_p$ random function is meant to be the
    It$\hat{o}$ integral. This integral is believed to be extended
    in the standard way to an isometric operator mapping
    $\mathfrak{L}_2 = L^2([0,T]\times \Omega, \overline{\mathcal{P}},\lambda_1 \times \P, \R)$
    into $L^2(\Omega, \mathcal{F}_T, \P, \R)$.  For an arbitrary
    function( equivalence class) in $\mathfrak{L}_2$ the value of
    the integral is , by definition, an equivalence class in $L^2(\Omega, \mathcal{F}_T, \P,
    \R$ containing the integral of a progressively measurable
    representative which always exists [2,p.11] in a class
    $\mathfrak{L}_2$. The stochastic integral in (2.3) is defined
    for every $t$ for a.e. $x$ as an element of $L^2(\Omega, \mathcal{F}_T, \P,
    \R)$.

\par
    THEOREM 2.1.  For any function $g\in \mathcal{H}^l$ a pair of
    functions $v$,$\mathcal{X}$, where $v \in X^1 \cap C_0 \cap
    C_2^{r+2}$, $r=[l]<l$, $\mathcal {X}=\|\mathcal {X}_1, \dots, \mathcal
    {X}_d\|$, $\mathcal{X}_j \in X^0$, $j=1,\dots,d$, is defined
    satisfying (2.1)-(2.2). Moreover, relation (2.2) holds for $t=T$
    for a.e.$(x,\omega)\in D \times \Omega$, and for $D\neq R^n$
    and $x\in \p D$ for a.e. $(t,\omega)\in [0,T] \times \Omega$.
    These functions $v0$, $\mathcal{X}_j$ are determined uniquely up
    to equivalence (as elements of $X^0$).

\par
    Let us note that the Bismut backward equations [5], which occur
    in the control theory for ordinary It$\hat{o}$ equations, have a
    form analogous to (2.1)-(2.2): one must find a solution of an
    It$\hat{o}$ equation adapted to a nondecreasing (unlike the backward equations of
    [2,p.36]) filtration of $\sigma$-algebras which takes on a given
    (for example, nonrandom) value at a finite time. Usually this
    problem is solvable for the only possible diffusion coefficient
    which must be found in the course of the solution (thus under the conditions of Theorem
    2.1 in view of uniqueness of $\mathcal{X}$ for nonrandom $f$, $\beta$, $g$, we have $\mathcal{X}\equiv
    0$). It$\hat{o}$ equations in an infinite-dimensional phase
    space, in particular parabolic It$\hat{o}$ equations, are by now
    well investigated (see, for example, [2] and [6]-[18] and their
    bibliographies). The corresponding infinite-dimensional
    analogues of the Bismut equations have practically not been
    studied at though they were introduced in [19].

\par
    We introduce the operators $\mathcal{T}$, $\mathcal {G}$,$\mathcal
    {G}_j$, $j=1,\dots,d$, which map a function $g$ into the
    functions $v=\mathcal{T}g$, $\mathcal{X}=\mathcal {G}g$, $\mathcal {X}_j = \mathcal {G}_j
    g$ , respectively, satisfying (2.1)-(2.2).

\par
    THEOREM 2.2  The operator $\mathcal {T}$ can be extended from
    the set $\mathcal {H}^l$ which is everywhere dense in $X^0$ and
    in $X^{-1}$ to continuous linear operators $\mathcal {T}: X^{-1}\rightarrow
    X^1$, $\mathcal {T}: X^{-1}\rightarrow C_0$, $\mathcal {T}: X^0 \rightarrow \mathcal
    {W}^2$.

\par
    In what follows the continuity of some operator signifies the
    possibility of its continuous extension from some everywhere
    dense set. The operators $\mathcal {T}$, $\mathcal {G}_j$ and
    others are the corresponding continuous extensions to $X^{-1}$
    (or in stipulated cases to $X^0$ or $X^{-1}$). An assertion of
    the type $"v = \mathcal {T}g\in C^0_2$ for $g\in \mathcal {H}^l$
    and operator $\mathcal {T}: X^{-1}\rightarrow X^1"$ means that
    $v$ and $g$ are representatives with the required properties for
    the functions (classes) $v = \mathcal{T}g \in X^1$, $g \in
    X^{-1}$.

\par
    THEOREM 2.3. The operators $\mathcal {G}: X^{-1}\rightarrow
    X^0$, $j=1, \dots , d$, are continuous.

\par
    DEFINITION. A generalized solution of the problem (2.1)-(2.2)
    for $g \in X^{-1}$ is a pair of functions $v$, $\mathcal{X}$,
    where $v=\mathcal{T}g\in X^1 \cap C_0$, $\mathcal {X}=\mathcal {G}g=\|\mathcal {X}_1,\dots,\mathcal
    {X}_d\|$, $\mathcal {X}_j \in X^0$.

\par
    THEOREM 2.4. The operators $\mathcal {G}_j: X^0 \rightarrow \mathcal
    {W}^1$, $j= 1,\dots, d$, are continuous. Moreover, $\mathcal {G}_j g \in
    X^1$ for $g \in X^0$ (we recall that $X^1 = \mathcal {W}^1 $ for $D = R^n$ and $X^1 \subset \mathcal{W}^1$).

\par
    THEOREM 2.5. Assume $f\in C^2_2$, $\beta \in C^2_2$, let the
    function $g\in C_2^2$ be a representative of some class in
    $X^{-1}$, and for the equivalence classes $\mathcal {T} g \in
    X^0$, $\mathcal {G}_j g \in X^0$, let there exist
    representatives $v \in C_2^4 \cap C_0$, $\mathcal {X}_j \in
    C_2^2$. Then the function
\[
   \phi(x,t,\omega) = g(x,t,\omega)-\sum_{j=1}^{d}(\beta_j(x,t,\omega) \frac {\p \mathcal{X}}{\p x}(x,t,\omega)
\]
    belongs to $\mathcal{C}_0^2$, and, for $v$ and $\phi$, relation
    (1.3) holds for any $s\in[0,T]$ for a.e. $(x,\omega)\in[0,T\times
    \Omega]$ for any $x\in D$.
    We introduce the operator $B: X^0 \rightarrow X^0$ by the
    formula
\begin{equation} \tag{2.4}
    Bg = -\sum_{j=1}^d \beta_j(x,t,\omega)\frac{\p \mathcal{X}_j}{\p
    x}(x,t,\omega), \text{where} \mathcal{X}_j = \mathcal{G}_j g
\end{equation}

THEOREM 2.6. Let $r>(n/2)+2$, $g\in \mathcal {H}^l$(the numbers $l$,
$r=[l]$ are fixed in Section 1). Then the hypotheses of Theorem 2.5
hold and there exist representatives $v\in \mathcal {C}^4_2 \cap
\mathcal{C}_2^0 $, $\varphi \in \mathcal {C}_2^0$ of the equivalence
classes $\mathcal{T}g \in X^0$, $g+Bg \in X^0$ for which relation
(1.3) holds for any $s$ for a.e.$x$, $\omega$ and for a.e. $s$,
$\omega$ for any $x$.
\par
    Thus the solution $\mathcal {T}g$ of the boundary value problem
    is represented in the form (1.3), where
    \begin{equation} \tag{2.5}
        \varphi = g+Bg
    \end{equation}
    The question arises whether(2.5) is solvable with respect to $g$
    for a given $\varphi$.
    \par
  \emph{  Proof of Theorem 2.1. }For a solution of the problem (1.5)-(1.6)
    we have $U\in \overline{\mathfrak{C}}^{r+2,1} \subset C([0,T])\rightarrow
    L^p(\Omega,\mathcal{F}_T, P, C^{r+2}(\overline{D}))$, $p=1,2$. For $\xi \in L^1(\Omega,
    \mathcal{F}_T, P, C^{r+2}(\overline{D}))$ the symbol
    $\xi_{\mathcal{F}}$. $\xi$ denotes the projection of $\xi[20]$
    to the space $L^1(\Omega, \mathcal{F}_s), P, C^{r+2}(\overline{D})$.
    We introduce the functions $v(x,t,\omega)=\xi_{\mathcal{F}_T} U(x,t,\omega)$
    and $u(x,t,s,\omega)=\xi_{\mathcal{F}_T} U(x,t,\omega)$. We have
    $v \in \mathcal{C}^{r+2}_1$, $u(.,s,.)\in
    \overline{\mathcal{C}}_1^{r+2}(s)$.
\par
    Below let the symbol $\mathcal{D}_x^1$ denote any partial
    derivative in $x$ of order $l$, $0\leq l \leq r+2$, and let the
    symbol $\mathcal{D}$ denote either $\mathcal{D}^l_x$ or $\p / \p
    t$.
\par
    By the Clark theorem(see [21,p.178]) we have the representation
    \begin{equation} \tag{2.6}
        \mathcal{D}U(x,t,\omega)=E \mathcal{D} U(x,t,\omega)+
        \sum_{j=1}^d \int_0^T r^{\mathcal{D}}_j(x,t,\rho,\omega) d
        \omega_j(\rho) \text{ a.s.}
    \end{equation}
\par
    Here $r^{\mathcal{D}}_j$ are some functions of the class
    $C(\overline{Q}\rightarrow \mathfrak{L}_2)$ (since
    $\mathcal{D}U \in C( \overline{Q}\rightarrow L^2(\Omega, \mathcal{F}_T,
    P,R))$); the order of arguments in such that
    $\overline{Q}=\{(x,t)\}$.
\par
    Let $\gamma_i$ denote the function in (2.6) for $\mathcal{D}U = U$
    (that is, $\gamma_i = \gamma_i^{\mathcal{D}_x^0}$). It can be
    easily sen that all other $\gamma_j^{\mathcal{D}}$ are the
    derivatives of the form $\mathcal{D}_{\gamma}$ of the functions
    $\gamma_j: \overline{Q}\rightarrow \mathfrak{L}_2$, and
    $\gamma_j \in C^{2,1}(\overline{Q} \rightarrow \mathfrak{L}_2)$.
    Below the partial derivatives $\mathcal{D}_{\gamma_j}$ which
    occur, for example, i the expression $A(x,t,\omega)\gamma_j(x,t,\rho, \omega)$
    are assumed to be the functions $r_j^{\mathcal{D}}$.
\par
    Let us prove that the function $v$ introduced above and the
    functions
    \begin{align*}
        \mathcal {X}_j(x,t,\omega)=r_j(x,0,t,\omega)- \int_0^t
        A(x,\rho,\omega)r_j(x,\rho,t,\omega)d\rho,\\
        \mathcal{X} = \|\mathcal{X}_1,\dots,\mathcal{X}_d\|,  \tag{2.7}
    \end{align*}
    are the ones required.

\par
    We have
    \[
        \mathcal{D}^l_xv(x,t,\omega) = E\{ \mathcal{D}_x^lU(x,t,\omega)|
        \mathcal{F}_t\}, \mathcal{D}u(x,t,s,\omega)=\E \{\mathcal{D}U(x,t,\omega)
        |\mathcal{F}_s\}
    \] a.s.

    $U\in X^l \cap \mathcal{C}_2^{r+2}(T)$. Thus the functions $v$
    and $\mathcal{D}_x^1 v$ are square-summable in $x$, $t$,
    $\omega$ and the function $\|v\|_{C^{r+2}(\overline{D})}$ in $t$,
    $\omega$, that is, $v \in X^1 \cap \mathcal{C}_2^{r+2}$.
    Obviously $v(x,t,\omega)-v(x,s,\omega)=\zeta_1 + \zeta_2$, where
    $\zeta_1 = u(x,t,t,\omega)-u(x,t,s,\omega)$ and
     $\zeta_2 = u(x,t,s,\omega) -u(x,s,s,\omega)$. Letting $t-s$
     tend to $0+$, we have
    \begin{align*}
        \E \|\zeta_1\|^2_{L_2(D)} \leq \E
        \|U(x,t,\omega)-U(x,s,\omega)\|^2_{L_2(D)} \rightarrow 0,\\
         \E \|\zeta_2\|^2_{L_2(D)} \leq \sum_{j=1}^d \E \|\int_s^t \gamma_j(x,s,
         \rho,\omega)^2d\rho\|_{L_1(D)} \rightarrow 0.
    \end{align*}
    Consequently, $v\in C_0$. By (1.6), relation (2.2) holds for
    a.e. $x$, $\omega$ for $t=T$ and for a.e. $t$, $\omega$ for $x\in \p
    D$, $D\neq R^n$. In (2.7) the coefficients of one derivatives in
    $A(x,\rho,\omega)$ are bounded, continuous, and
    $\mathcal{F}_{\rho}-$adapted for a.e.$\omega$; hence $\mathcal{X}_j \in
    X^0$.
\par
    By virtue of (1.5)-(1.6) and (2.6) we have, for a.e.
    $x$,$\omega$,
    \begin{align*}
        v(x,t,\omega)   = &\E\{U(x,t,\omega)| \mathcal{F}_t\}\\
                        = &u(x,0,t,\omega)-\int_0^t[A(x,s,\omega)u(x,s,t,\omega)+g(x,s,\omega)]
                            ds\\
                        = &v(x,0,\omega) + \sum_{j=1}^d \int_0^t
                        \gamma_j(x,0,\rho,\omega) d\omega_j(\rho)\\
                          & - \int_0^t \{A(x,s,\omega)[v(x,s,\omega)+ \sum_{j=1}^d \int_s^t \gamma_j(x,s,\rho,
                          \omega)d\omega_j(\rho)]+g(x,s,\omega)\}ds\\
                        = &v(x,0,\omega)- \int_0^t[A(x,s,\omega)v(x,s,\omega)+
                        g(x,s,\omega)]ds\\
                        &+ \sum_{j=1}^d\{\int_0^t \gamma_j(x,0,\rho,\omega)d\omega_i(\rho)
                        -\int_0^t ds \int_s^t A(x,s,\omega)\gamma_j(x,s,\rho,\omega)
                        d\omega_j(\rho)\}.
    \end{align*}
    The sum of the It\^{o} integrals in the right-hand side of the
    latter equality is equal to $\int_0^t
    \mathcal{X}(\rho)d\omega(\rho)$ by (2.7) and the Fubini theorem
    for stochastic integrals (see[22]). So, for $v$, $\mathcal{X}$
    relation (2.3) holds and $v$, $\mathcal{X}$ are the ones
    required.
\par
    We introduce the operator $\mathcal{T}^{*}:X^0\rightarrow X^0$
    by the rule $\mathcal{T}^*h=\pi$, where the function $\pi \in X^1 \cap C_0$
    is a solution of the boundary value problem
    \begin{equation} \tag{2.8}
        \frac {\p \pi}{\p t}(x,t,\omega) =
        A^*(x,t,\omega)\pi(x,t,\omega) + h(x,t,\omega),
        \pi(x,t,\omega)|_{(x,t)\in \p_0 Q}=0.
    \end{equation}

    The operator $\mathcal{T}: X^0 \rightarrow X^0$(and even the operator $\mathcal{F}^*:X^{-1}
    \rightarrow X^1$) is linear and continuous (see [2] and [3]).
    The dual operator in the Hilbert space $X^0$ is denoted by
    $\mathcal{T}$; the operator $\mathcal{T}: X^0 \rightarrow X^0$
    is continuous. For some $v^{'} \in C^2_2 \cap C_0 \cap X^1$ and
    $\mathcal{X}^{'}_j \in X^0$, let (2.1)-(2.2) hold as indicated
    in the theorem. It can be verified immediately that
    $(\mathcal{T}^*h,g)_{X^0}=(h,h^{'})_{X^0}(\forall h \in X^0)$.
    So $v^{'}=\mathcal{T}g$ in $X^0$ and hence $v=v^{'}$ in $X^0$.
    From (2.3) we obtain that, if $v, v^{'} \in C^2_2\cap C_0$ and
    $v^{'}=v$ in $X^0$, then $\mathcal{X}_j=\mathcal{X}^{'}_j(\forall j)$
    in $X^0$. Thus $v$ and $\mathcal{X}_1,\dots,\mathcal{X}_d$ are
    determined uniquely in $X^0$. The theorem has been proved.
\par
    The proof of Theorem 2.2 follows from the estimates
    \[
        \|v\|_{X^1}+ \|v\|_{C^0}\leq \|U\|_{\overline{X}^1} + \|U\|_{C^0}
        \leq c_1 \|g\|_{X^1}, \|v\|_{\omega^2} \leq c_2 \|g\|_{X^0}
    \]
    which hold for constants $c_i > 0$ common for all $g$, $v$, $U$
    in the proof of Theorem 2.1, by virtue of known (see [2]-[4])
    properties of the operators $\overline{\mathcal{T}}$ and properties
    of the operation $\E\{.|\mathcal{F}_t\}$.
\par
    The proof of Theorem 2.3 will be adduced in Section 4.
\par
    \emph{Proof of Theorem 2.4. } Let us consider functions
    \[
        \gamma = \gamma(x,t,\rho,\omega)\in L^2([0,T]\times \Omega,
        \overline{\mathcal{P}}, \lambda_1 \times \P, C^2(\overline{Q})) \cap
        L^2([0,T]\times \Omega, \overline{\mathcal{P}}, \lambda_1 \times \P, W^2_2(Q))
    \],
    which are equal to zero in the case $D\neq \R^n$ for $x \in \p
    D$ for a.e. $\rho$, $\omega$(here $Q = \{(x,t)\}$).
    Using the estimate [4,p.523,(149)] for $D \neq \R^n$ and a
    similar estimate for $D=\R^n$, for a constant $c_1 >0$ common
    for all such $\gamma$, we obtain the estimate
    \begin{align*}
        \E \int_0&T \|\gamma(x,t,t,\omega)\|^2_{W^1_2(D)}dt \leq \E
        \int_0^T \sup_{t\in
        [0,T]}\|\gamma(x,t,\rho,\omega)\|^2_{W61_2(D)}d \rho \\
        &\leq c_1 \E \int_0^T(\|\frac{\p \gamma}{\p t}(x,t,\rho, \omega)\|^2 _{L_2(Q)}+ \sum_{i=1}^n
        \|\frac{\p \gamma}{\p x_i}(x,t,\rho,\omega)\|^2_{L_2(Q)}\\
        &+ \sum_{t,j=1}^n \|\frac{\p^2 \gamma}{\p x_i \p x_j}\|^2_{L_2(Q)}
        )d\rho. \tag{2.9}
    \end{align*}

\par
    Obviously this estimate can be extended to all functions $\gamma = \gamma(x,t,\rho, \omega)$
which belong to $C^{2,1}(\overline{Q} \rightarrow \mathfrak{L}_2)$,
are square-summable in $\overline{Q}\times [0,T]\times \Omega $
together with the corresponding derivatives and are equal to zero
for $x\in \p D$ intthe case $D \neq \R^n$. Such are the functions
$\gamma_i$ in the representation (2.6) for $U= \mathcal{T}g$, $g\in
\mathcal {H}^l$. The right-hand side of the latter inequality in
(2.9) under the substitution $\gamma = \gamma_j$ is finite and dose
not exceed the value
    \begin{equation}
     c_2\left(\|U\|^2_{\overline{\mathcal{W}}^2}+ \|\frac{\p U}{\p t}\|^2_{\overline{X}^0}\right)
        \leq c_3 \|g\|^2_{X^0}, \tag{2.10}
    \end{equation}

where $c_i >0$ are constants the same for all $g\in \mathcal{H}^l$.
Thus (2.11)
    \begin{equation} \tag{2.11}
        \|\gamma_j(x,t,t,\omega)\|_{\mathcal{W}^1} \leq
        \sqrt{c_3}\|g\|_{X^0}.
    \end{equation}

    From (1.5)-(1.6) and (2.6) we obtain
    \begin{align*}
        \sum_{j=1}^d \int_s^T \frac{\p \gamma_j}{\p
        s}(x,s,\rho,\omega)d \omega_j(\rho) & = \frac{\p U}{\p
        s}(x,s,\omega)- \E \{\frac{\p U}{\p s}|\mathcal{F}_s\}\\
            &=
            -A(x,s,\omega)[U(x,s,\omega)-\E\{U(x,s,\omega)|\mathcal{F}_s\}]\\
            &= -\sum_{j=1}^d \int_s^T A(x,s,\omega)
            \gamma_j(x,s,\rho,\omega)d\omega_j(\rho).
    \end{align*}

    for $(x,s)\in Q$ with probability 1. This relation and (2.7)
    imply that, for a.e. $x,t,\omega$,
    \[
        \mathcal{X}_j(x,t,\omega)= \gamma_j(x,0,t,\omega) + \int_0^t
        \frac {\p \gamma_j}{\p s}(x,s,t,\omega)ds =
        \gamma_j(x,t,t,\omega).
    \]
    By extending the estimate (2.11) from the everywhere dense
    subset $\mathcal{H}^l$ to $X^0$, we obtain the assertion of the
    theorem.
\par
   \emph{ Proof of Theorem 2.5.} As is seen from (2.1), the differential
    \[
        d_t v(x,t,\omega) = \tilde{v}(x,t,\omega)dt +
        \mathcal{X}(x,t,\omega)d\omega(t)
    \],
    where $\tilde{v}=-Av-g \in C62_2$,$\mathcal{X}\in C^2_2$,exists.
    We assume that $\tilde{v}(x,t,\omega)$ and $\mathcal{X}(x,t,\omega)$
    are defined on $\R^{n+1} \times \Omega$ and equal to zero for
    $(x,t,\omega)\notin \overline{Q} \times \Omega$. For $\epsilon >0$ we
    introduce the functions
    \begin{align*}
        \tilde{v}_{\epsilon}(x,t,\omega) &= \epsilon ^{-1} \int_{t-\epsilon}
        ^t \tilde{v}(x,\rho,\omega)d\rho,\\
        \mathcal{X}_{\epsilon}(x,t,\omega) &= \epsilon^{-1} \int_{t-\epsilon}^t \mathcal{X}(x,\rho,\omega)d\rho,
        \mathcal{X}_{\epsilon}=\|\mathcal{X}^{(\epsilon)}_1,\dots,
        \mathcal{X}_d^{(\epsilon)}\|,\\
        v_{\epsilon}(x,t,\omega)&= \E U(x,0,\omega)+
        \int_0^t\tilde{v}_{\epsilon}(x,\rho,\omega)d\rho + \int_0^t
        \mathcal{X}_{\epsilon}(x,\rho,\omega)d \omega(\rho).\\
    \end{align*}
    We denote $\gamma^{x,s}(\epsilon, \omega)=(T+\epsilon)\wedge \inf
    \{t: y^{x,s}(t,\omega) \notin \overline{D}\}$. From the
    It$\hat{o}$-Venttsel formula(see [2] and [23]), whose
    applicability  is left without a proof, we see that for a
    modification of the function $v_{\epsilon}$ for $(x,s) \in Q$ in
    the class $\overline{\mathfrak{C}}^{2,0}$ the following relation
    holds:
    \[
        v_{\epsilon}(x,s,\omega) = -\E\left\{ \int_s^{r^{x,s}(\epsilon, \omega)}
        \left[ \tilde{v}_{\epsilon} + A\tilde{v}_{\epsilon} + \sum_{j=1}^d\beta_j
        \frac{\p \mathcal{X}_j^{(\epsilon)}}{\p x} \right]
        (y^{x,s}(t,\omega),t,\omega)dt|\mathcal{F}_s\right\} a.s.
    \]
    Restricting all functions again to $\overline{Q} \times \Omega$, we
    have, as $\epsilon \rightarrow 0$,
    \[
        \tilde{v}_{\epsilon} \rightarrow -Av-g = \tilde{v},
        Av_{\epsilon} \rightarrow Av, v_{\epsilon} \rightarrow v,
        \beta_j \frac{\p \mathcal{X}_j^{(\epsilon)}}{\p x}
        \rightarrow \beta_j \frac{\p \mathcal{X}_j}{\p x}
    \]
    in the metric of $C^0_2, v_{\epsilon} \rightarrow v$ in the
    metric of $C_0$. In addition $\gamma^{x,s}(x,\omega)\rightarrow
    \gamma ^{x,s}(\omega)$ uniformly in $\omega \in \Omega$. Hence
    we obtain the assertion of the theorem.

\par
    \emph{Proof of Theorem 2.6.} We have $v \in C^{r+2}_2\cap C_0$, $v(x,0,\omega)=\E U(x,0,\omega)$
    for a.e. $x$, $\omega$,$\E U(x,0,\omega) \in C^{r+2}(\overline{D})$.
    From this and also from (2.3), with $t=0$, as well as from the
    Clark theorem([21,p.178]) we obtain successively for $l=0,1,\dots, r-1$
    for arbitrary $i,j$ and the vector $e_i =
    \|\delta_{k_i}\|_{k-1}^n$(where $\delta_{k_i}$ is the Kronecker
    symbol)that the limit of the expression
    \[
        \epsilon ^{-1}\{\mathcal{D}^l_x\mathcal{X}_j(x+\epsilon e_i,t,\omega)-
        \mathcal{D}^l_x \mathcal{X}_j(x,t,\omega)\},
    \]
    as $\epsilon \rightarrow 0$,exists in $X^0$ which we denote from
    now on by $\p \mathcal{D}^l_{x}\mathcal{X}_j(x,t,\omega)$. Hence
    we can approximate $\mathcal{X}_j(x,t,\omega)$ by functions
    $\mathcal{X}_j^{(i)}\in \mathcal{W}^r$ so that
    $\mathcal{D}^L_x \mathcal{X}^{(i)}_j \rightarrow \mathcal{D}^l_x \mathcal{X}_j$
    in $X^0$ as $i\rightarrow +\infty$, $l = 0,1,\dots,r$(we can use averagings
    in $x$ of the type of [1,p.48] with a smooth kernel for a.e. $t$,$\omega$,
    extending $\mathcal{X}_j$ to $\R^n \times[0,T]\times \Omega$ for $D \neq \R^n$).
    The completeness of $W^r$ implies the existence of the limit $\mathcal{X}_j^{(i)}$
    in $W^r$ equivalent to $\mathcal{X}_j$ in $X^0$. From the
    inclusion (see [3,p.61]) $\mathcal{W}\subset C^2_2$ we obtain
    the  required assertion.

\section{Forms and Properties of dual operations}

    In addition to the operation $\mathcal{T}$, $\mathcal{G}_j$, $B$
    introduced above we shall consider the operators $R=
    (I+B)^{-1}$, $L=\mathcal{T}R$. The operator $R$ maps a function
    $\varphi$ into a solution $g=R\varphi$ of the equation (2.5)
    connected with the problem (2.1)-(2.2); $I$ is the identity
    operator. The operator $L$ maps the function $\varphi$ into the
    joint solution $v=L \varphi=\mathcal{T}(R\varphi)$ of the
    equation (2.5) and the problem(2.1)-(2.2).

\par
    The symbols $\mathcal{T}^{*}$,
    $\mathcal{G}_j^{*}$,$B^{*}$,$R^{*}$,$L^{*}$, and so on denote
    the corresponding dual operators in the Hilbert space $X^0$ (we shall show that
    the operators $R$ and $L$ are well defined on sets which are everywhere dense in  $X^0$)
\par
    For an n-vector $\xi=\|\xi_i\|_{i=1}^n$ we denote $(\nabla,\xi)=\sum_{i=1}^n \p \xi_i/ \p x_i$.
\par
    Below we shall consider initial-boundary value problems of the
    type of [2,$\S\S$3.4-4.1] with a boundary condition at $t=0$.
    The symbol $\mathcal{X}^0$ will denote the set of processes
    $h(x,t,\omega)$ which are representatives of some functions in
    $X^0$, predictable [2,p.16] for all $x$, and taking values in
    $L_2(D)$ for all $t$,$\omega$. The symbol $\mathcal{X}^{-1}$ will denote the set of  and
    processes $h$which are representatives of
    functions (classes) in $X^{-1}$ and representable in the form
    $h=(\nabla,\xi)$, where $\xi=\|\xi_1,\dots,\xi_n\|$, $\xi_i \in \mathcal{X}^0(\forall i)$ .
    Solutions of boundary value problems are defines in [2] for free
    terms in $\mathcal{X}^k$. It is known [24, Chap.3] that in every
    equivalence class of $X^{-1}$, $X^0$ there are representatives
    of $\mathcal{X}^{-1}$,$\mathcal{X}^0$, respectively. Therefore,
    we can (and shall) understand by a solution of boundary value
    problems of the type of [2] with an initial condition at $t=0$
    for free terms in $X^k$ ,$k=-1,0,$ and extension to these
    Hilbert spaces of continuous operators (using suitable theorems of
    [2]) which map free terms of boundary value problems into
    solution in $X^1 \cap C_0$. Then a boundary condition of the
    form $g(x,t,\omega)|_{(x,t)\in \p_0 Q}=0$ is said to be
    satisfied if $g\in X^1 \cap C_0$ and $g(x,0,\omega)=0$ for a.e.
    $x,\omega$.

\par
    THEOREM 3.1. The operators $\mathcal{G}_j^{*}: X_0 \rightarrow
    X_1$ are continuous and have the form $\mathcal{G}_j^{*}h = q$,
    where the function $q\in X^1 \cap C_0$ satisfies the boundary
    value problem
    \begin{equation}\tag{3.1}
        d_t q(x,t,\omega)= A^{*}(x,t,\omega)q(x,t,\omega)dt +
        h(x,t,\omega)d \omega_j(t),
    \end{equation}

    \begin{equation}\tag{3.2}
        q(x,t,\omega)|_{(x,t)\in \p_0 Q}=0.
    \end{equation}
\par
    THEOREM 3.2. The operator $B^{*}:X^0 \rightarrow X^0$ is
    continuous and has the form $B^{*}h = z$, where the function $z$
    satisfies the boundary value problem

    \begin{equation}\tag{3.3}
        d_t z(x,t,\omega) = A^{*}(x,t,\omega)dt +
        \sum_{j=1}^d(\nabla,\beta_j(x,t,\omega)h(x,t,\omega))dw_j(t),
    \end{equation}

    \begin{equation}\tag{3.4}
        z(x,t,\omega)|_{(x,t)\in \p_0 Q} = 0.
    \end{equation}

    For $h \in \mathcal{X}^1$ the solution $z = B^{*}h \in X^1 \cap
    C_0$ is understood in the sense of [2], for $h \in X^0$, $h \notin
    \mathcal{X}^1$, the solution is the limit in $X^0$ of a sequence
    $B^{*}h_i$, where $h_i \in \mathcal{X}^1$ and $\|h_i-h\|_{X^0} \rightarrow 0$
    as $i \rightarrow + \infty$.
\par
    The theorem stated above contains the assertion of existence of
    a "generalized" solution in the class $X^0$(or of the possibility of defining a solution as the corresponding limit in this
    space) for a coefficient belonging to the class $X^{-1}$ of the
    stochastic differential in the free term of the equation. This
    assertion is apparently new for the theory of partial
    It$\hat{o}$ equations.
\par
    THEOREM 3.3. For $d<d_0$, the operator $R^{*}: X^0 \rightarrow
    X^0$ is determined uniquely and the operator $R^{*}: X^1 \rightarrow X^1 $
    is continuous. For $\pi \in X^1$, the operator has the form
    $R^{*}\pi = h$, where $h= \pi -z$ and the function $z \in X^1 \cap
    C_0$ is a solution of the boundary value problem
    (3.5)
    \begin{align*}
        d_t z (x,t,\omega) =& A^{*}(x,t,\omega)z(x,t,\omega)dt \\
                                & +
                                \sum_{j=1}^d(\nabla,\beta_j(x,t,\omega)[\pi(x,t,\omega)-z(x,t,\omega)])d\omega_j(t),
                                \tag{3.5}
    \end{align*}

    \begin{equation} \tag{3.6}
        z(x,t,\omega)|_{(x,t)\in \p _0 Q}=0.
    \end{equation}
\par
    THEOREM 3.4. For $d<d_0$, the operator $L^{*}:X^1 \rightarrow
    X^1$ is continuous and has the form $L^*\xi = h$, where the
    function $h \in X^1 \cap C_0$ is a solution of the boundary
    value problem

    \begin{align*}
        d_t h(x,t,\omega) = &[A^*(x,t,\omega)h(x,t,\omega)+ \xi
        (x,t,\omega)]dt \\
            & - \sum_{j=1}^d(\nabla,
            \beta_j(x,t,\omega)h(x,t,\omega))d \omega_j(t), \tag{3.7}
    \end{align*}

    \begin{equation}\tag{3.8}
        h(x,t,\omega)|_{(x,t)\in \p_0Q} =0.
    \end{equation}
\par
    Let us note that (3.1) and (3.3) are superparabolic [2]
    It$\hat{o}$ equations, and (3.5)and (3.7) are superparabolic for
    $d<d_0$ and parabolic for $d=d_0$.
\par
    Proof of Theorem 3.1. First let $f$ and $\beta$ be nonrandom.
\par
    Suppose that $g\in \mathcal{H}^l$ is an arbitrary function and
    the functions $\hat{g}_j \in C(\overline{Q} \rightarrow \mathfrak{L}_2)$
    are determined by the Clark theorem [21,p. 178] from the
    representation

    \begin{equation}\tag{3.9}
        g(x,t,\omega) = \E g(x,t,\omega) + \sum_{j=1}^d \int_0^t
        \hat{g}_j(x,t,\rho,\omega) d \omega_j(\rho);
    \end{equation}
    the functions $u(x,t,s,\omega)\in C^2_2(s)$ and $\gamma_j(x,t,\rho,\omega)\in C^{2,1}(\overline{Q}\rightarrow \mathfrak{L}_2)$
    for $U = \overline{\mathcal{T}}g \in \overline{\mathfrak{C}}^{2,1}$ are
    defined in the same way as in the proof of Theorem 2.1. We have
    \be\ba
        \sum_{j=1}^d \int_0^t \frac{\p \gamma_j}{\p
        t}(x,t,\rho,\omega)d \omega_j(\rho) = \E \{\frac{\p U}{\p t}(x,t,\omega)|
        \mathcal{F}_t\} - \E \frac{\p U}{\p t}(x,t,\omega)\\
         = -\left[A(x,t)u(x,t,t,\omega)+ \E\{ g(x,t,\omega)| \mathcal{F}_t\}-A(x,t)u(x,t,0,\omega)-\E g(x,t,\omega)\right]
   \nonumber \\
        = - \sum_{j=1}^d \int_0^t [A(x,t)\gamma_j(x,t,\rho,\omega)+
        \hat{g}_j(x,t,\rho,\omega)]d \omega_j(\rho).
    \ea\eqno{(3.10)}\ee
    Let $G(x,y,t,s)$ be Green's function of the boundary value
    problem (1.5), (1.6) with the nonrandom operator $A(x,t,\omega)=
    A(x,t)$; then from (3.10) and the condition $\gamma_j(x,t,\rho,\omega)|_{(x,t)\in \p_0 Q}=0$
    for a.e. $\rho$, $\omega$ we obtain that, for a.e. $\rho$,
    $\omega$,

    \[
        \gamma_j(x,t,\rho,\omega) = \int_D dy \int_t^T
        G(x,y,t,s)\hat{g}_j(y,s,\rho,\omega)ds.
    \]

\par
    For an arbitrary function $h \in \mathcal{H}^l$ we have
    \begin{align*}
    (\mathcal{G}_j g, h)_{X^0}&= \E \int_0^T dt \int _D dx
    \left\{\left[\int_D dy \int_t^T G(x,y,t,s)\hat{g}_j(y,s,t,\omega)ds\right]h(x,s,\omega)\right\}
    \\
    &= \E \int_0^T ds \int_Ddy \int_0^s dt
    \hat{g}_j(y,s,t,\omega)\int_D G(x,y,t,s)ds \\
    &= (g,\mathcal{G}^*_jh)_{X^0}.
    \end{align*}
    In view of (3.9) and the fact that $g$ is arbitrary this means
    that
    \[
        (\mathcal{G}^*_j h)(y,s,\omega)= \int_0^s d\omega_j(t)
        \int_D G(x,y,t,s)h(x,t,\omega)dx.
    \]
    From this relation we obtain (3.1), (3.2) and the form of $\mathcal{G}^*_j$
    for nonrandom $f$,$\beta$.
\par
    Now let $f, \beta$ be random. Consider the function $f_0(x,t)= \E f(x,t,\omega)$
    ,$\beta_0(x,t)= \E\beta(x,t,\omega)$. Let $A_0, A^*_0, \mathcal{T}_0, \mathcal{G^*_{j,0}}$
    denote the operators corresponding to the function $f_0,
    \beta_0$, which are defined like the $A,A^*,
    \mathcal{T},\mathcal{G}_j$ are defined for the functions $f,
    \beta$. We introduce the operator $\mathfrak{A} = (A_0-
    A)\mathcal{T}_0$: for $g_0 \in X^0$ we have $\mathfrak{A} g_0 =
    (A_0-A)v_0$, where $v_0=\mathcal{T}_0 g_0$. The operators
    $\mathfrak{A}: X^{-1} \rightarrow X^1$ and $\mathfrak{A}: X^0 \rightarrow X^0$
    are continuous by Theorem 2.2.

\par
    It can be verified immediately that, for $g= g_0 +
    \mathfrak{A}g_0$, we have $\mathcal{T}g = \mathcal{T}_0 g_0$ and
    $\mathcal{G}_j g = \mathcal{G}_{j,0}g_0$. This means that
    $\mathcal{G}_j = \mathcal{G}_{j,0}(I +\mathfrak{A})^{-1}$ and
    the dual operator in $X^0$ has the form $\mathcal{G^*}_j = (I+\mathfrak{A})^{-1}
    \mathcal{G}^*_{j,0}$.

\par
    Obviously $\mathfrak{A}^* = \mathcal{T}^*_0(A^*_0 - A^*)$. The
    form of $\mathcal{T}^*$ (and analogously of $\mathcal{T}^*_0$)
    was established in $\S2$ by formulas (2.8). The operators
    $\mathcal{T}^*$ and $\mathcal{T}^*_0$ map $X^{-1}$ continuously
    (see [2]-[4]) into $X^1$ and $C_0$. The operators $A^*$ and
    $A^*_0$ map $X^{1}$ continuously into $X^{-1}$.

\par
    For $q \in X^1$ we have $z = \mathfrak{A}^* q \in X^1 \cap C_0$,
    and $z= z(x,t,\omega)$ satisfies the boundary value problem
    $d z/ dt = A^*_0 z + (A^*_0- A^*)q, z|_{(x,t)\in \p _0Q = 0}$ in
    Q.
\par
    Let us find the form of $a = (I +\mathfrak{A}^{-1})\eta$ for $\eta \in
    X^1$. We have $q+\mathfrak{A}q = \eta$. We denote $z= \mathfrak{A}^*
    q$; then $q = \eta -z$ and $z= \mathfrak{A}^*(\eta-z)$. By the
    assertion proved above the function $z=z(x,t,\omega)$ satisfies
    the boundary value problem
    \begin{equation} \tag{3.11}
        \frac{dz}{dt}=A^*_0 z + (A_0^*-A^*)(\eta-z) = A^* z + (A^*_0 -
        A^*)\eta, z|_{(x,t)\in \p_0 Q} = 0.
    \end{equation}

    Thus, $z=\eta-q \in X^1 \cap C_0$ and $q \in X^1$.
\par
    Let us find the form of $q = (I + \mathfrak{A}^*)^{-1} \eta$ for
    $\eta = \mathcal{G}^*_{j,0} h$, where $h \in X^0$. By Theorems
    3.4.8 and 4.1.1 of [2] the function $\eta = \eta(x,t,\omega)\in X^1 \cap C_0$
    and satisfies, by the proof, the boundary value problem
    $d_t \eta = A^*_0 \eta dt + h d\omega_j(t), \eta|_{((x,t)\in \p_0
    Q}=0$. Moreover, $q= \eta -z $ where the function $z=
    z(x,t,\omega)$ satisfies a boundary value problem of the form
    (3.11). So  $q \in X^1 \cap C_0$. From the formulas for $d_t
    \eta$, and $d_t z = (dz/dt)dt$ we find $d_t q = d_t \eta -d_t z$
    and thus we obtain (3.1). Condition (3.2) is satisfied since the
    analogous conditions hold for $z$ and $\eta$. Continuity of the
    operator $\mathcal{G}^*_j: X^0 \rightarrow X^1$ (and even continuity of the operator
    $\mathcal{G}^*_j: X^0 \rightarrow C_0$) follows from Theorem
    3.4.8 and 4.1.1 of [2]. Theorem has been proved.
\par
    \emph{Proof of Theorem 3.2.} By Theorem 2.4, for $g \in X^0$, $\mathcal{X}_j g \in X^1$
    . For $h \in X^1$ we have
    \[
        (Bg, h)_{X^0} = \sum_{j=1}^d (\mathcal{X}_j, (\nabla, \beta_j
        h))_{X^0} = \left(g, \sum_{j=1}^d \mathcal{G}^*_j(\nabla, \beta_j h )\right)_{X^0}
    \]
    From this relation and linearity of the problem (3.1), (3.2) we
    obtain (3.3),(3.4). Theorem 3.4.8 and 4.1.1 of [2] imply
    continuity of the operators $B^*: X61 \rightarrow X^1$ and
    $B^*: x^1 \rightarrow C_0$. Continuity of the operator $B: X^0 \rightarrow
    X^0$ proved in Theorem 2.4 implies continuity of the operator $B^*: X^0 \rightarrow
    X^0$.

\par
    \emph{Proof of Theorem 3.3.} If $h = R^* \pi$, then $h = \pi -z$, where
    $z = B^*h = B^*(\pi - z)$. By substituting the value $h = \pi -
    z$ into (3.3)-(3.4) we obtain formulas (3.5), (3.6). continuity
    of the operator $R^*: X^1 \rightarrow X^1$ (and thus uniqueness of the operator $R^*:
    X^0 \rightarrow X^0 $ ) follows from Theorem 3.4.8 and 4.1.1 of
    [2].
\par
   \emph{ Proof of Theorem 3.4.} Let $h = R^*\pi$ and $z = B^*h$; then $h=
    \pi - z$. Let $\pi = \mathcal[T]^* \xi$ be determined from the
    problem (2.8), where $\xi \in X^{-1}$. Using (2.8) and (3.3),
    (3.4) we obtain the expression for $d_th = d_t\pi - d_tz$ or
    (3.7), as well as condition (3.8). Continuity of the operator
    $L^* : X^{-1} \rightarrow X^1$ follows from the form of $h =
    L^*\xi$ and Theorems 3.4.8 and 4.1.1 of [2].

\section{Solvability of (2.5)}

    THEOREM 4.1. Let $d< d_0$. Then the operator $R: X^0 \rightarrow X^0$
    is unique and well-defined on some everywhere dense set $\mathcal{D}(R) \subset X^0$
    in $X^0$( and in $X^{-1}$) (that is, (2.5) has at most one solution $g \in X^0$
    for any $\phi \in X^0$ and , in addition, the set of those $\phi \in X^0$, for which the equation
    is solvable with respect to $g \in X^0$, is everywhere dense in
    $X^0$). The operators $R$ and $L $ defined on $\mathcal{D}(R)$
    can be extended from this set to operators defined on $X^{-1}$
    so that the operators $R: X^{-1} \rightarrow X^{-1}$, $L: X^{-1}\rightarrow X^1$
    , and $L:X^{-1}\rightarrow C_0$ are continuous.
\par
    Let $\mathcal{C}_*$ denote the set of all $\phi \in C^0_2 \cap
    X^0$ such that up to equivalence (in $X^0$) $\phi = g + Bg$ for
    some $g \in \mathcal{H}^l$, when $l >r$, the integer number $ r \geq
    0$, and the number $l$ were introduced in $\S1$. We recall that
    such $\phi$ occurred in the statement of Theorem 2.6, which
    asserts that for these $\phi$ with $r>n/2+2$ there exists a
    modification in the class $C^0_2$(and hence in the class
    $\mathcal{C}_*$) and moreover for this modification the value of
    the functional (1.3) coincides with $L \phi$.
\par
    THEOREM 4.2.  For $d> d_0$ and $r> n/2 + 2$, the set
    $\mathcal{C}_*$ is everywhere dense in $X^0$.
\par
    \emph{Proof of Theorem 4.1.} The operator $R^*$ is defined on the set
    $X^1$ which is everywhere dense in $X^0$ and maps $X^1$. Hence
    the operator $I+B^*$ inverse to it has a set of values
    everywhere dense in $X^0$ and  a kernel consisting only of zero.
    Obviously, the operator $I+B$ has the same properties. Thus the
    operator $R$ is determined uniquely and has a domain that is
    everywhere dense in $X^0$ and is denoted by $\mathcal{D}(R)$.
    Since $X^0$ is everywhere dense in $X^{-1}$, $\mathcal{D}(R)$ is
    everywhere dense in $X^{-1}$.
\par
     It remains to prove the assertion concerning the extension of
     the operators to $X^{-1}$. For $k=0, \pm1, \pm2, u = u(x,t,\omega)\in X^0$
     the symbol $\Lambda^k u $ will denote the function in $X^{-k}$
     obtained by the application, for a.e. $t, \omega$ of the
     operator $\Lambda^k$ introduced in $\S 1$ in the definition of
     the spaces $H^k$ to the function $u(.,t,\omega)\in L_2(D)$. In
     view of Theorem 3.3 we have, for some constant $c>0$ and for
     any $\varphi \in \mathcal{D}(R), h \in X^0$,
     \par
    \begin{align*}
    (R\varphi, h)_{X^{-1}} &= (\Lambda^{-1} R \varphi, \Lambda^{-1} h)_{X^0} =
    (R\varphi, \Lambda^{-2} h)_{X^0} = ( \Lambda^{-1}\varphi, \Lambda R^*
    \Lambda^{-2}h)_{X^0}\\
    &\leq \|\varphi\|_{X^1}\|R^*\Lambda^{-2}h\|_{X^{-1}}\leq c
    \|\varphi\|_{X^{-1}}\|\Lambda^{-2}h\|_{X^1} \leq
    c\|\varphi\|_{X^{-1}}\|h\|_{X^{-1}}
    \end{align*}
This inequality implies the existence of a continuous extension of
the operator $R$ to $X^{-1}$. The corresponding assertion of the
theorem for the operator $L = \mathcal{T}R$ follows from Theorem
2.2. The theorem has been proved.

\par
    \emph{Proof of Theorem 4.2.} For an arbitrary number $\xi >0$ and for
    $\varphi \in X^0$, it is required to find $\varphi_0 \in \mathcal{C}_*$
    such that $\|\varphi - \varphi_0\|_{X^0}< \xi$. By theorem 4.1
    there exists $\varphi^{'} \in X^0$ such that $\varphi^{'} = g^{'} + B g ^{'}$
    for some $g^{'} \in X^0$ and $\| \varphi- \varphi^{'}\|< \xi/2$.
    The norm $\|I+B\|$ of the operator $(I+B): X^0\rightarrow X^0$
    is positive since this operator has an image which is everywhere
    dense in $X^0$. For $g^{'}\in R\varphi$, there exists $g^{"}
    \in \mathcal{H}^l$ such that $\|g^{'}-g^{"}\|_{X^0}<
    \xi\|I+B\|^{-1}/2$. By Theorem 2.6, the function $\varphi_0 = g^{"} +
    Bg^{"}$ has a modification in the class $\mathcal{C}^0_2$. This
    function $\varphi \in \mathcal{C}^0_2$ is the one required since
    \[
        \|\varphi - \varphi_0\|_{X^0} \leq \|\varphi -
        \varphi^{'}\|_{X^0}+ \|\varphi^{'}-\varphi_0\|_{X^0} < \xi/2
        + \|g^{'}-g^{"}\|_{X^0}\|I+B\|< \xi
    \]
    The theorem has been proved.
    \par
    \emph{Proof of Theorem 2.3.} By virtue of Theorem 3.1, we have, for
    some constant $c>0$ and for any $g \in X^0$, $h \in X^0$,
    \begin{align*}
        (\mathcal{G}_j g , h)_{X^0} & = (g, \mathcal{G}^*_j h) \leq
        \|\Lambda ^{-1} g\|_{X^0}\|\Lambda \mathcal{G}^*_j
        h\|_{X^0}\\
        & = \|g\|_{X^{-1}}\|\mathcal{G}^*_j h\|_{X^1} \leq
        c\|g\|_{X^{-1}} \|h\|_{X^0}
    \end{align*}.
    This inequality implies the assertion of the theorem.

    \section {Representation of functionals of It$\hat{o}$ processes in
    the form of solutions of boundary value problems}

    Let us adduce some sufficient conditions for the functional
    (1.3) to coincide, for a given $\varphi$, with a solution of the
    problem (2.1), (2.2) and (2.5).

    \par
    THEOREM 5.1. Let $D= \R ^n $ or $D\neq \R^n, d< d_0$; let the
    function $\beta(x,t,\omega)= \beta(t,\omega)$ not depend on $x$,
    the function $f \in \C_2^2$, and let at least one of the
    following conditions hold:

    \par
    a) the function $f(x,t,\omega)=f(t,\omega)$ does not depend on
    $x$ and the function $\varphi \in \C_2^0 \cap X^0$;
    \par
    b) $n = 1$ and the function $\varphi \in \C^0_2 \cap X^0$;
    \par
    c) $n =1$ and $\varphi = \varphi(x,t)$ is a nonrandom Borel
    measurable function of $L_2(Q)$.

    \par
    Then the value $v(x,s,\omega)$ of the functional (1.3) as a
    function of $(x,s,\o)$ belongs to $X^1 \cap C_0$ and coincides
    with $L \varphi$ as a function in $X^0$ and in $C_0$(i.e., is a solution
    of the problem (2.1),(2.2), and (2.5)).

    \par
    COROLLARY 5.1. Under the assumptions of Theorems 4.1 and 5.1 the
    estimates
    \[
        \|v\|_{C_0} \leq c_1 \|\varphi\|_{X^{-1}} \leq c_2 \|
        \varphi\|_{X^0}, \|v\|_{X^1} \leq c_3
        \|\varphi\|_{X^{-1}}\leq c_4 \|\varphi\|_{X^0},
    \]
    hold for the functional (1.3), where the constants $c_i>0$
    depend only on $n,d,d_0, Q$, and the values

    \begin{align*}
        \delta = \inf_{x,t,\o} Det
        \tilde{\beta}(x,t,\o)\tilde{\beta}(x,t,\o)^T, K_1 =
        \sup_{x,t,\o} |f(x,t,\o)|,\\
        K_2 = \sup_{x,t,\o}|\beta(x,t,\o)|, K_3^i =
        \sup_{x,t,\o}\left|\frac{\p \beta}{\p x_i}(x,t,\o)\right|
    \end{align*}
    (cf. estimates in [1,  $\S\S\Pi.2-\Pi.3$]).

    \par
    \emph{Proof of Theorem 5.1.} Let assumptions a) or b) hold. For a
    function $\eta(x) \in L_2(\R^n)$ the symbol $(\eta)_{\xi}$ will
    denote its averaging (convolution) with the kernel of the
    Sobolev averaging $\zeta(x/\xi)\xi^{-n}$. Here the function $\zeta(x)=0$
    for $|x|\geq 1$, $\zeta(x)= \mathcal{X}_n \exp \{|x|^2(|x|^2-1)\}$
    for $|x|<1$; $\mathcal{X}_n $ is a normalizing factor such thar
    $\int_{\R^n}\zeta(x)dx =1$/. If $\eta \in H^{-1}$ and $\eta = \p \xi / \p x_j$
    , where $\xi = L_2(\R^n)$, $j \in \{1,\dots, n\}$, then we
    assume that
    \[
        (\eta)_{\xi}(x) = - \xi ^{-n-1} \int_{\R^n}\frac{\p \zeta}{\p x_j}
        \left(\frac{y-x}{\xi}\right)\xi(y)dy.
    \]
\par
    Let $D= \R^n$. For functions $u \in X^0$ or $u \in X^{-1}$, the
    symbol $(u)_{\xi}$ denotes a function in $C^0_2$ coinciding
    with $(u(., t,\o))_{\xi}$ for all$t,\o$ such that $u(.,t,\o) \in H^0 = L_2(\R^n)$
    or $u(.,t,\o) \in H^{-1}$, respectively (i.e.,for a.e.,$t,\o$).
\par
    Let $D \neq \R^n$. In this case, functions defined on $\overline{Q}\times
    \Omega$ are assumed to be extended to $\R^n \times [0,T]\times
    \Omega$, and the operation $(.)_{\xi}$ is applied to them
    according to the rule indicated above.

\par
    Everywhere in the proof of this theorem, $C^2_m$ will be the
    space $C^2_m = L_2([0,T]\times \Omega), \overline{\mathcal{P}}, \lambda_1 \times P, C^m(\R^n)$
    . So, for $m = 0,1,2,\dots$ and $u \in X^0$ or $u \in X^{-1}$,
    we have $(u)_{\xi} \in C^2_m$ for the spaces $X^0, X^{-1}$
    defined for $D= \R^n$ as well as $D \neq \R^n$.
\par
    For $D \neq \R^n$ we denote by $D_{\xi}$ a region with a
    $C_2$-smooth boundary which contains the union of
    $2\xi$-neighborhoods for all $x\in D$ and which itself is
    contained in the union of $3\xi$-neighborhoods for all $x\in D$
    . The symbol $\tau_{\xi}^{x,s}(\omega)$ denotes the random time
    $T\wedge \inf\{t: y^{x,s}(t,\o)\notin \overline{D}_{\xi}\}$. For $D= \R^n$
    we assume that $D_{\xi}=D=\R^n$,
     $\tau_{\xi}^{x,s}(\o)\equiv\tau^{x,s}(\o)\equiv T$.
\par
    Let $D=\R^n$ or $D \neq \R^n$, $g = \R \varphi \in X^{-1}$,
    $v=L\varphi$, $\mathcal{X}=\mathcal{G}g$, $\mathcal{X} _j =
    \mathcal{G} _j g$. Then $v \in X^1 \cap C_0, \mathcal{X}_j \in
    X^0$ and, for all $s, x$, we have
\begin{align*}
    (v)_{\xi}(x,s,\o) = \int_s^T(Av+g)_{\xi}(x,t,\o)dt
    -\int_s^T(\mathcal{X})_{\xi}(x,t,\o) d \o(t),\\
    (v)_{\xi}(x,t,\o)|_{x\in \p D_{\xi}}= 0 \text{  in the case } D \neq \R^n
    , v_{\xi}(x,T,\o)=0
\end{align*}
    with probability 1.

\par
    We introduce the functions $\Delta_{\xi}= (Av)_{\xi} - A(v)_{\xi}, \Phi_{\xi }= \Delta
    _{\xi}+ (\varphi)_{\xi}$. These functions belong to the class
    $C^2_2$. The function $(v)_{\xi} \in C^4_2$ is a solution of the
    problem of the form (2.1), (2.2) with the free term $(g)_{\xi} + \Delta_{\xi}
    \in C^2_2 \cap X^{-1}$ in the cylinder $D_{\xi}\times (0,T)$. By
    Theorem 2.5, for any $s \in [0.T]$ for a.e.$(s,\o) \in [0,T] \times
    \Omega$, we have
\begin{align*}
    (v)_{\xi}(x,s,\o) = & \E \left\{ \int_s^{\tau_{\xi}^{x,s}(\o)} \Phi_{\xi}[y^{x,s}(t,\o), t, \o]dt
    |\mathcal{F}_s\right\}\\
                      = & \E\left\{\int_{\tau^{x,s}(\o)}^{\tau_{\xi}^{x,s}(\o)} \Phi_{\xi}[y^{x,s}(t,\o), t, \o]dt |\mathcal{F}_s\right\}\\
                        & + \E \left\{ \int_s^{\tau_{\xi}^{x,s}(\o)} \Phi_{\xi}[y^{x,s}(t,\o), t, \o]dt
    |\mathcal{F}_s\right\} \tag{5.1}
\end{align*}
    for any $x \in D_{\xi}$.

\par
    Let us estimate $\Delta_{\xi}$. For a.e. $t, \o$ we have $v(.,t,\o)\in
    H^1$ and
    \[
        \Delta_{\xi}(x,t,\o) = \xi^{-n} \int_D \frac{\p v}{\p
        x}(y,t,\o)[f(y,t,\o)- f(x,t,\o)] \zeta\left(\frac{x-y}{\xi}\right)dy.
    \]
    Under assumption a) this value is equal to zero.

\par
    The H\"older inequality and the inequality $\zeta(x)^2 < \mathcal{X}_n \zeta(x)$
    imply that

\begin{align*}
    &|\Delta_{\xi}(x,t,\o)|\\
    &\leq c_1 \xi^{-n} \left(\sup_{|x-y|\leq \xi}|f(x,t,
    \o)-f(y,t,\o)|\right)\|v(.,t,\o)\|_{H^1} \left(\int_{R^n}\zeta(\frac{x-y}{\xi})^2
    dy\right)^{1/2}\\
    & \leq c_1 \xi ^{-n} c_2 \xi \left( \sup_{x,t,\o}|\frac{\p f}{\p x}(x,t,\o)|\right)
    \|v(.,t,\o)\|_{\mathcal{H}^1} \xi ^{n/2} = \xi_1(t, \o) \xi ^{n/2-n+1}
\end{align*}

    for a.e.$t, \o$ where $c_i >0$ are constants and $\xi_1(t,\o)$
    is some function in $L^2([0,T]\times \Omega), \overline{\mathcal{P}}, \lambda_1 \times P, \R$
    .
\par
    Let $\xi \rightarrow 0$. For the function
    \[
        \xi_2(t,\o) = 2 \sup_{x\in D}|\varphi(x,t,\o)| \in L^2([0,T]\times \Omega,
        \overline{\mathcal{P}}
        ,\lambda_{1} \times \P, \R),
    \]
    we have $|(\varphi)_{\xi}(x,t,\o)|+ |\varphi(x,t,\o)| \leq
    \xi_2(t,\o)$ for a.e. $t,\o$. Moreover, for a.e. $t,\o$,
    $(\varphi)_{\xi} (x,t,\o)\rightarrow \varphi(x,t,\o)$ for any $x \in
    D$.
    We have also $\Delta_{\xi} \rightarrow 0$ in the metric of
    $C^0_2$ for $n = 1$. Thus,
    \begin{align*}
    &\E \int_0^T \left|(\varphi)_{\xi}[y^{x,s}(t,\o),t, \o] - \varphi[y^{x,s}
    (t,\o ),t,\o]\right|^2 dt \rightarrow 0,\\
    &\E \int_0^T \left|\Delta_{\xi}[y^{x,s}(t,\o), t, \o]\right|^2 dt
    \rightarrow 0.
    \end{align*}
From these relations and the Lebesgue theorem we see the first term
in the right-hand side of equality (5.1) tends, in the metric of
$L^2(\Omega, \mathcal{F}, \P, \R)$, to the right-hand side of
equality (1.3). Moreover, the left-hand side of (5.1) tends to $v= L
\varphi$ in the metrics of $C_0$ and $X^0$ as a function of $x,s,\o$

\par
    To complete the proof of the theorem for the case of assumptions
    a) and b) we prove that the second term in the right-hand side
    of (5.1) tends to zero in the metric of $L^2(\Omega, \mathcal{F}, \P,
    \R)$. Obviously this term is equal to zero in case $D=\R^n$. For
    $D \neq \R^n$ and $n=1$, we obtain
    \[
        \Phi_{\xi}[y^{x,s}(t,\o), t, \o] \leq \xi_1(t, \o) \xi^{1/2}
        + \xi_2(t, \o)
    \]
    for a.e.$t,\o$. Moreover, $\tau_{\xi}^{x,s}(\o)\downarrow
    \tau^{x,s}(\o)$ a.s, since $\tau^{x,s}(\o) = T \wedge \inf\{
    t: y^{x,s}(t,\o) \notin \overline{D}\}$ and for a.e.$\o$ there exists
    $\overline{\theta} = \overline{\theta}(\o) >0$ such that $y^{x,s}(t, \o) \in \R^n \ \overline{D
    }$ for $\tau^{x,s}(\o) < T$, $t = \tau ^{x,s}(\o) +
    \theta, \theta \in ( 0, \overline{\theta}(\o) ]$. Hense we obtain the
    required assertion for assumptions a) and b).

\par
    Let assumption c) hold. We introduce the operator $\tilde{L}$
    defined on $L_2(Q)$, mapping functions $\varphi \in L_2(Q)$ into
    the values $\tilde{v}=\tilde{L}\varphi$ of the functional(1.3)
    regarded as functions of $(x,s, \o)$.
\par
    Assume $\varphi \in L_2(Q)$, $\tilde{v}=\tilde{L}\varphi$, the
    sequence $\{\varphi_i\}_{i=1}^{+\infty} \subset C(\overline{Q})$,
    $\tilde{v}_i=\tilde{L}\varphi_i$, and $\varphi_i \rightarrow \varphi$
    in the metric of $L_2(Q)$ as $i \rightarrow + \infty$. By the
    above proof, $\tilde{L}\varphi_i = L \varphi_i \in C_0$. Theorem
    $ \Pi.2.4$ and $\Pi.3.4$ of [1] imply that, for some constant
    $c>0$,
    \[
        \sup_{t \in [0,T]}\E \|\tilde{v}(x,t,\o)-v_i(x,t,
        \o)\|^2_{H^0} \leq c\|\varphi-\varphi_i\|_{L_2(Q)}.
    \]
    Completeness of the space $C_0$ implies that $\tilde{v}_i \rightarrow
    v$ in $C_0$ as $t \rightarrow +\infty$. Hence $\tilde{v} = \tilde{L}\varphi = L \varphi$

\par
    Theorem 5.1 has been proved.

\section{On distributions of It\^o processes}

    Let $D = \R^n $ or $D \neq \R^n$, $d<d_0$ and let $p_0(x) \in
    L_2(D)$ be some nonrandom functions. We consider in the cylinder
    $Q$ the boundary value problem

    \begin{equation}\tag{6.1}
        d_i p(x, t, \o) = A^*(x,t, \o) p (x,t ,\o) dt - \sum_{j=1}^d
        (\nabla, \beta_j(x,t, \o) p(x,t, \o))d\o_j(t),
    \end{equation}

    \begin{equation}\tag{6.2}
        p(x,0, \o) = P_0(x), p(x,t,\o)|_{x\in \p D = 0}.
    \end{equation}

    The boundary condition on $\p D$ in (6.2) is not considered for
    $D = \R ^n$.
\par
    Equation (6.1) is a superparabolic It\^o equation[2]. A
    solution of the problem (6.1)-(6.2) is understood to be
    analogous to [2]; this problem has a solution $p \in X^1 \cap
    C_0$. The boundary conditions (6.2) for $p \in X^1 \cap C_0$ are
    said to be satisfied if $p(x,0, \o)= p_0(x)$ for a.e. $x, \o$.
\par
    LEMMA 6.1. For $\phi \in X^0$ and $s \in [0,T]$, the equality

    \begin{equation} \tag{6.3}
        \int_D p(x,s, \o) v(x,s,\o)dx = \E \left\{ \int_s^T dt \int_D p(x,t,\o)dx |\mathcal{F}_s\right\}
    \end{equation}
    holds with probability 1. Here $v= L\varphi \in X^1 \cap C_0$ is
    a solution of the problem (2.1), (2.2), (2.5).
\par
    Let, in (1.1), (1.2), $s=0$ and $x = x (\o)$ be a random
    n-vector. We assume that $\E|x(\o)|^2 < + \infty$, $x(\o)\in D$
    a.s. the vector $x(\o)$ does not depend on $W(t)$ for any $t \geq
    0$ and has a probability density $p_0(x) \in L_2(D)$. Let $y^{x(\omega), 0}(t,\o)$
    be the corresponding solution of equations (1.1),(1.2) and let
    the random time $\tau^{x(\o), 0}(t,\o)= T \wedge \inf \{t: y^{x(\o),0}\} \notin \overline{D}$
    The symbol $I_{\tau}(t,\o)$ denotes the indicator function of
    the event $\{ \tau ^{x(\omega), 0} \leq t\}$. For $D = \R^n$ we
    have $\tau^{x(\o), 0}(\o)(\o)= T$, $I_{\tau}= (\tau, \o)\equiv
    1$, for $0 \leq t \leq T$.

\par
    THEOREM 6.1. Let $\varphi \in \overline{\mathfrak{C}}\cap C^0_2 \cap
    X^0$, and let the assumptions of Theorem 5.1 hold. Then for a.e.
    $(t, \o)\in [0,T] \times \Omega$( and even for any $t\in [0,T]$ almost surely
    if $D = \R^n$) the following equality holds:

    \begin{equation}\tag{6.4}
        \E \{I_{\tau}(t, \o)\varphi[y^{x(\o), 0}(t,\o), t,
        \o]|\mathcal{F}_t\} = \int_D p(x,t, \o)\varphi(x,t, \o)dx.
    \end{equation}
\par
    This theorem establishes the distribution of the process $y^{x(\o),0}(t, \o)$
    (broken off at the exit of $\overline{D}$ if $D\neq \R^n$); $\E p (x,t,\o)
    $ is the distribution density of the process. A close result is
    proved in [2, Th. 5.3.1], where equality (6.4) is obtained for
    $D = \R^n$ and coefficients $f, \beta$ of general form (no restrictive assumptions of
    Theorem 5.1 are required). Moreover, in [2] another method of
    the proof is used, and equality (6.4) is obtained only for
    nonrandom $\varphi$ and $D=\R^n$, which is essential. Theorem
    5.3.1 of [2] establishes the distribution of the It\^o process
    $y^{x(\o), 0}(t,\o)$; therefore with its help we can obtain the
    following analogue of Theorem 5.1(less strong, however, for functions $f$ and
    $\beta$ of general form).

\par
    THEOREM 6.2. Let $D= \R^n, d<d_0$, and let the function
    $\varphi(x,t)\in C(\overline{Q})\cap L_2(Q)$ be nonrandom. Then
    \begin{equation}\tag{6.5}
        \E \int_0^T \varphi[y^{x(\o),0}(t,\o), t]dt = \int_{\R^n}
        p_0(x) \overline{v}(x, 0)dx = \E \overline{v}[x(\o), 0].
    \end{equation}
    Here $\overline{v}(x,0) \in L_2(\R^n)$ is a nonrandom modification of
    the function $v(x,s,\o)|_{s=0}$ where $v = L \varphi \in X^1 \cap C_0$
    ( in other words, $v(x, 0)= v(x, 0, \o)$ for a.e. $x, \o$).
\par
     \emph{The proof of Lemma 6.1 }follows from equalities (2.4), (2.5) and
     the equality
     \begin{align*}
     \E &\{(p(x,s,\o), v(x,s,\o))_{H^0}|\mathcal{F}_s\}\\
      = & \E \{ (p(x, T, \o), v(x, T, \o))_{H^0} - \int_s^T dt [(A^*(x,t, \o)p(x,t,\o),
      v(x,t,\o))_{H^0}     \\
      & - (p(x,t,\o), A(x,t,\o)v(x,t,\o)+ \sum_{j=1 }^d \frac{\p \mathcal{X}_j}{\p x
      }(x,t,\o)\beta_j(x,t,\o)+ \varphi(x,t,\o))_{H^0}\\
      & - \sum_{j=1}^d ((\nabla, \beta_j(x,t,\o)p(x,t,\o)),\mathcal{X}_j(x,t,\o))_{H^0}]| \mathcal{F}_s    \}.
     \end{align*}

\par
   \emph{ Proof of Theorem 6.1.} Let $\xi \in L^2([0,\tau]\times \Omega, \overline{
    \mathcal{P}}, \lambda \times \P, \R)$ be an arbitrary function .
    Consider the functions $\tilde{\varphi}(x,t,\o)= \varphi(x,t,\o)\xi(t,\o)$
    and $v= v(x,s, \o)= L \tilde{\varphi} \in X^1 \cap C_0$. We have
    $v(x,0, \o)\in L^2(\Omega, \mathcal{F}_0, \P, L_2(D))$. The
    probability of any event of $\mathcal{F}_0$ is equal to 0 or 1.
    Thus the function $v(x, 0, \o)$ has a nonrandom modification $\tilde{v}(x,0)
    \in L_2(D)$ such that $v(x, 0) = v(x, 0, \o)$ for a.e. $x, \o$.

\par
    For $x \in D$ we consider the $(n+1)$-dimensional process
    \[
        \eta^x (t, \o) = \| y^{x,0}(t, \o), z^x(t, \o)\|, \text{ where
        } z^x (t, \o) = \int_0^t \tilde{\varphi}[y^{x,0}(\rho, \omega), \rho, \omega]
        d \rho.
    \]
    Analogously we define the process $\eta^{x(\omega)}(t, \o)$ for
    a random vector $x(\o)$ using the process $y^{x(\o), 0}(t, \o)$
    instead of $y^{x, 0}(t, \omega)$.
\par
    On functions of the form $\eta(t) \|y(t), z(t)\|$, where $y(.) \in C([0,\tau]
    \rightarrow \R^n)$ and $z(.) \in C([0,\tau]\rightarrow \R)$, we
    define the functional $F[\eta(.)]= z(\tau)$, where $\tau = \min \{T,
    \inf\{t: y(t)\notin \overline{D}\}\}$. By Theorem 5.1, $\tilde{v}(x,0) = \E f[\eta^x(., \o)]$
    for a.e. x. By virtue of Theorem $\Pi.9.4$ of [1] establishing
    an analogue of the Markov property for It\^o processes, we have
    $\E \tilde{v}[x(\o), 0] = \E F[\eta^{x(\o)}(.,\o)]$. Thus

    \[
        \int_D p_0(x)\tilde{v}(x,0)dx = \E \tilde{v}[x(\o), 0] = \E
        \int_0^{\tau^{x(\o),0}(\o)} \tilde{\varphi}[y^{x(\o), 0}(t, \o),t
        ,\o] \xi (t, \o) dt.
    \]

    From these equalities and equality (6.3), where $s=0$, we obtain
    \begin{align*}
    \E &\int_0^T \left( \int_D p(x,t, \o)\varphi(x,t,\o)dx\right)\xi
    (t, \o) dt\\
    &= \E \int_0^T I_{\tau}(t, \o) \varphi[y^{x(\o), 0}, t, \o]
    \xi(t, \o)dt.
    \end{align*}
    Since $\xi$ is arbitrary, this relation implies the assertion of
    the theorem.

\par
    \emph{Proof of Theorem 6.2.} The existence of a nonrandom modification
    for the function $v(x,0, \o)$ can be established as in the proof
    of Theorem 6.1. By virtue of Theorem 5.3.1 of [2] the left-hand
    member of equality (6.5) coincides with $\E\int_Q p(x,t,\o) \times \varphi(x,t)
    dxdt$. By Lemma 6.1 this value is equal to the middle member of
    equality (6.5) (and hence to the right-hand member of this
    equality). The theorem has been proved.

\subsection*{Acknowledgment} The author thanks I. A. Ibragimov
for his interest in this paper.
\subsection*{References}

$\hphantom{xk}$
[1] N. V. KRYLOV, Controlled Diffusion Processes,
Springer-Verlag, New York, 1980.
\par
[2] B. L. ROZOVSKII, Evolutionary Stochastic Systems. Linear
Theory and Applications to Statistics of Stochastic Processes,
Nauka, Moscow, 1983. (In Russian.)
\par
[3] O. A. LADYZENSKAJA, V. A. SOLONNIKOV, AND N. N. URALTSCEVA,
Linear and Quasi-Linear Equations of Parabolic Type, American
Mathematical Society, Providence, RI, 1968.
\par
[4] V. I. SMIRNOV, Course of higher mathematics, Vol. IV, part 2,
Nauka, Moscow, 1981. (In Russian.)
\par
[5] J. M. BISMUT, Conjugate convex functions in optimal stochastic
control, J. Math. Anal. Appl., 44 (1973), pp. 384-404.

\par
[6] V. V. BAKLAN, On existence o] solutions of stochastic equations
in a Hilbert space, Dopovodi AN URSR, 10 (1963), pp. 1299-1303. (In
Ukrainian.)
\par
[7] YA. I. BELOPOL'SKAYA AND YU. L. DALETSKII, Diffusion processes
in smooth Banach spaces and manifolds, I. Trudy Mosc. Math. Ob-va,
37 (197S), pp. 107-141. (In Russian.)
\par
[8] M. I. VISHIK, A. I. KOMECH, AND A. V. FURSIKOV, Some
mathematical problems of statistical hydromechanics, Russ. Math.
Surveys, 34 (1979), pp. 149-234.
\par
[9] I. I. GIKHMAN, A boundary value problem for a stochastic
equation of parabolic type, Ukrainian Math. J., 31 (1979), pp.
483-489. (In Russian.)
\par
[10] I. I. (IKHMAN, On a mixed problem for a stochastic
differential equation of parabolic type, Ukrainian Math. J., 32
(1980), pp. 367-372. (In Russian.)
\par
[11] N. V. KRYLOV AND B. L. ROZOVSKIJ, Stochastic evolution
equations, J. Sov. Math., 16 (1981), pp. 1233-1277.
\par
[12] L. G. MAIGULIS, Nonlinear filtering of bounded diffusion
processes and boundary value problems for partial stochastic
differential equations, Markov stochastic processes and their
applications, Vol. 1., Edition of Saratov University, Saratov, 1980,
pp. 50-63. (In Russian.)
\par
[13] S. YA. MAKHNO, Boundary value problems for partial stochastic
equations, Theory of stochastic processes, Vol. 12., Naukova Dumka,
Kiev, 1984, pp. 48-56. (In Russian.)
\par
[14] S. A. MELNIK, On smoothness of solutions of stochastic
differential equations of parabolic type, Theory of stochastic
processes, Vol. 11., Naukova Dumka, Kiev, 1983, pp. 78-82. (In
Russian.)
\par
[15] T. M. MESTECHKINA, Kolmogorov equation ]or solutions of
Cauchy problems for a class of linear evolution equations,
Ukrainian Math. J., 40 (1988), pp. 67-70.
\par
[16] A. V. SKOROKHOD, Random Linear Operators, D. Reidel, Dordrecht,
1984.
\par
[17] H. KUNITA, Cauchy problem for stochastic partial differential
equations arising in nonlinear filtering theory, System Control
Lett., 1 (1981), pp. 37-41.
\par
[18] E. PARDOUX, Stochastic partial differential equations and
filtering of diffusion processes, Stochastics, 3 (1979), pp.
127-167.
\par
[19] H. KWAKERNAAK, A minimum principle for stochastic control
problems with output feedback, System Control Lett., 1 (1981), pp.
74-77.
\par
[20] J. DIESTEL, Geometry of Banach Spaces. Selected Topics,
Springer-Verlag, Berlin, 1975.
\par
[21] I. I. GIKHMAN AND A. V. SKOROKHOD, Stochastic Differential
Equations and their Applications, Naukova Dumka, Kiev, 1982. (In
Russian.)
\par
[22] G. KALLIANPUR AND C. STRIEBEL, Stochastic differential
equations occurring in the estimation of continuous parameter
stochastic processes, Theory Probab. Appl., 14 (1969), pp. 567-594.

\par
[23] B. L. ROZOVSKII, On the Ito-Venttsel formula, Moscow Univ.
Math. Bull., 28 (1973), pp. 22-26.
\par
[24] K. L. CHUNG AND R. J. WILLIAMS, Introduction to Stochastic
Integration, Birkhuser, Boston, 1983.

\end{document}